\documentclass[12pt,a4paper,final]{amsart}

\usepackage{url} 
\usepackage[hidelinks]{hyperref}
\usepackage[latin1]{inputenc}
\usepackage{amssymb,amsthm,amsmath,amstext,amsfonts}
\usepackage{amsaddr}
\usepackage{graphicx}
\usepackage[a4paper,twoside=true,top=2.0cm,bottom=2.2cm,left=2.5cm,right=2.5cm]{geometry}
\usepackage{color}
\usepackage{lineno}
\usepackage[normalem]{ulem}

\def\du{\ensuremath{\mathrm{d}}}

\theoremstyle{definition}

\def\texttiny#1{{\text{\tiny{#1}}}}
\def\DC{{}^{\texttiny{C}}\! D}

\begin{document}

\title[Neglecting nonlocality leads to unreliable numerical methods] 
{Neglecting nonlocality leads to \\unreliable numerical methods for \\fractional differential equations}

%
\author{Roberto Garrappa}
\address[Roberto Garrappa]{Dipartimento di Matematica, Universit\`a degli Studi di Bari\\ Via E. Orabona 4 - 70126 Bari - Italy\\Member of the INdAM Research group GNCS}
\email[R.~Garrappa]{roberto.garrappa@uniba.it}
\thanks{This is the preprint of a paper published in [Commun. Nonlinear Sci. Numer. Simul., 70 (2019) 302-306] at \url{https://doi.org/10.1016/j.cnsns.2018.11.004}. The plots have been corrected according to Corrigendum published in [Commun. Nonlinear Sci. Numer. Simul.] at \url{https://doi.org/10.1016/j.cnsns.2019.02.021}. Work supported by the INdAM-GNCS under a 2018 project.}




\maketitle

\begin{abstract}
In the paper titled ``New numerical approach for fractional differential equations'' by A. Atangana and K.M. Owolabi [Math. Model. Nat. Phenom., 13(1), 2018], it is presented a method for the numerical solution of some fractional differential equations. The numerical approximation is obtained by using just local information and the scheme does not present a memory term; {{moreover, it is claimed that third-order convergence is surprisingly obtained by simply using linear polynomial approximations}}. In this note we show that methods of this kind are not reliable and lead to completely wrong results since the nonlocal nature of fractional differential operators cannot be neglected. We illustrate the main weaknesses in the derivation and analysis of the method in order to warn other researchers and scientist to overlook this and other methods devised on similar basis and avoid their use for the numerical simulation of fractional differential equations.
\end{abstract}

\section{Introduction}

In the paper titled ``New numerical approach for fractional differential equations'', authors A. Atangana and K.M. Owolabi and published in [Math. Model. Nat. Phenom., 13(1), 2018], it is proposed a numerical method (see \cite[Eq. (3.11)]{AtanganaOwolabi2018}) for solving fractional differential equations (FDEs) of (presumably) order $0< \alpha< 1$
\begin{equation}\label{eq:FDE}
	\DC^{\alpha}_0 y(t) = f(t,y(t))
\end{equation}
with $\DC^{\alpha}_0$ the fractional derivative of Caputo's type. This method presents two surprising features:
\begin{itemize}
	\item there is no memory term and the approximation of the solution is obtained by using just local information;
	\item the numerical approximation, which is obtained by means of linear interpolant polynomials, is claimed to converge to the exact solution with  order $3$ under the simple assumption that $f(t,y(t))$ is bounded. 
\end{itemize}

{{It is a well-known fact that fractional derivatives are nonlocal operators and, as discussed in some recent papers \cite{Tarasov2016,Tarasov2018}, nonlocality is an essential, distinctive and inviolable feature of fractional-order operators.}} 

{{Obtaining approximations without having to compute a persistent memory term of increasing length is an ambitious goal in order to save computation and techniques for incorporating just a part of the memory term, known under the name of \emph{short} or \emph{fixed memory principle}, have been investigated by several authors (e.g., see \cite{Deng2007,FordSimpson2001,Podlubny1999}). However, the naive truncation of the memory term to a fixed interval leads to a severe loss of accuracy \cite{FordSimpson2001} and an appropriate selection of the length of the fixed memory (which must be sufficiently long) is necessary to keep the error under control; reducing the memory to very local information (just the last two grid-points are involved in the method proposed in \cite{AtanganaOwolabi2018}) unavoidably leads to catastrophic and unreasonable results.}}


Furthermore, achieving the third order of convergence on the basis of linear interpolant polynomials would be an amazing property as well which, however, contradicts the basic knowledge in numerical analysis.

The claim of both features (local approximation and third-order convergence) is obviously made on the basis of weak reasonings and after illogical and wrong derivations. In this note we illustrate the main weaknesses in the derivation and analysis of the method in \cite{AtanganaOwolabi2018} and we show that methods of this kind are not reliable and lead to completely wrong results. 

The aim of this note is to warn other researchers and scientists to overlook such kind of methods, to prevent from their use for the numerical simulations of fractional differential equations and to avoid in the future similar errors in the derivation and in the analysis of the methods.

\section{Analysis of the method}

To help the reader to understand the way in which the method in \cite[Eq. (3.11)]{AtanganaOwolabi2018} is derived, we start by observing that it is obtained after considering the Volterra integral equation formulation of (\ref{eq:FDE})

\[
	y(t) = y(0) + \frac{1}{\Gamma(\alpha)} \int_0^t (t-\lambda)^{\alpha-1} f(\lambda, y(\lambda)) \du \lambda ,
\]
and subtracting the values of $y(t)$ at $t=t_n$ from the one at $t=t_{n+1}$ ($t_n = nh$ is the grid with constant step-size $h>0$), namely
\[
	y(t_{n+1}) - y(t_n) = \frac{1}{\Gamma(\alpha)} \int_0^{t_{n+1}} (t_{n+1}-\lambda)^{\alpha-1} f(\lambda, y(\lambda))  \du \lambda  - \frac{1}{\Gamma(\alpha)} \int_0^{t_n} (t_n-\lambda)^{\alpha-1} f(\lambda, y(\lambda))  \du \lambda  ,
\]
{{in order to represent the solution at time $t_{n+1}$ in terms of some increments of the same solution at the time $t_n$, an approach already explored for instance in \cite{Deng2007}.}}

{{Unlike what proposed in \cite{Deng2007}, where a shorter but sufficiently large memory is considered after an in-depth analysis, the authors of \cite{AtanganaOwolabi2018} replace the function $f$}} in each integral in the right-hand-side (rhs) of (\ref{eq:FDE}) by the linear interpolant polynomial on the nodes $t_{n-1}$ and $t_n$
\[
	f(\lambda,y(\lambda)) \approx f(t_{n},y_{n}) \frac{(\lambda-t_{n-1})}{h} - f(t_{n-1},y_{n-1}) \frac{(\lambda-t_{n})}{h} 
\]
and hence analytically evaluate the integrals to obtain the resulting method 
\[
	\begin{aligned}
	y_{n+1} =& \, \, y(t_n) + \frac{f(t_n,y_n)}{h \Gamma(\alpha)} \left\{ \frac{2h}{\alpha(\alpha+1)} t_{n+1}^{\alpha} -  \frac{t_{n-1} t_{n+1}^{\alpha}}{\alpha+1} - \frac{h}{\alpha}t_n^{\alpha} + \frac{t_n^{\alpha+1}}{\alpha+1} \right\}  \\
	&+ \frac{f(t_{n-1},y_{n-1})}{h \Gamma(\alpha)} \left\{ \frac{t_{n+1}^{\alpha+1}}{\alpha+1} - \frac{h}{\alpha} t_{n+1}^{\alpha} - \frac{t_n^{\alpha+1}}{\alpha+1} \right\}  \\
	\end{aligned}
\]
(this technique is known as product-integration (PI) rule \cite{Young1954} and we have corrected here a series of errors and misprints of the original text which will be highlighted later). 

A method of this kind is clearly not reasonable since the function in the rhs is approximated by the same polynomial on the whole integration interval $[0,t_n]$ or $[0,t_{n+1}]$. The integration interval in both integrals begins at 0 and the presence of weight functions $(t_{n}-\lambda)^{\alpha-1}$ or $(t_{n+1}-\lambda)^{\alpha-1}$ with real powers does not allow to recast the difference of the two integrals as an integral on the local interval $[t_n,t_{n+1}]$; as it is widely known, the remainder of interpolating polynomials depends on the distance from the nodes (with linear interpolation the remainder increases quadratically and in this case is proportional to $(t-t_{n})(t-t_{n-1})$), with $t$ varying, in this case, on the whole intervals $[0,t_n]$ and $[0,t_{n+1}]$. Therefore, as the integration moves forward (i.e. as $n$ increases) larger and larger errors are produced; building methods on the basis of approximations which deteriorate as the integration process proceeds clearly leads to poorly accurate methods.

This is the reason why in PI rules for integrals with singular kernels the integration interval is partitioned (usually by means of the same grid-points) and the polynomial approximation is built locally, in a piece-wise way, by using in each sub-interval $[t_{j-1},t_{j}]$ an interpolant polynomial on just nodes next to $t_{j-1}$ and $t_{j}$. Example of methods of this kind are described not only in \cite{Young1954} but also in several more recent papers (e.g. \cite{DiethelmFordFreed2002,DiethelmFordFreed2004,Dixon1985,Garrappa2015_MCS}) or in the survey \cite{Garrappa2018_Mathematics}.

It is therefore extremely surprising that, despite the poor and unreliable way by which the method is derived, the third order of convergence is proved in \cite[Theorem 3.1]{AtanganaOwolabi2018}, namely the difference at $t=t_n$ between the exact and approximated solution is claimed to be
\[
		R^{\alpha}_n < \frac{h^{3+\alpha} M}{12 \Gamma(\alpha+1)} \bigl\{ (n+1)^{\alpha} + n^{\alpha} \bigr\}
\]
which is a quantity proportional to $h^{3} t_{n+1}^{\alpha}$ (also in this case, to avoid confusion, we have corrected some misprints in the original statement of the Theorem).

Obviously, Theorem 3.1 is wrong and its proof is made after a series of inconsistencies and illogical derivations. In particular, for the benefit of readers, we observe that:
\begin{enumerate}
	\item In the second equation at page 6 the remainder $R_1(t)$ is defined in terms of $n$ but it is not clear in this context the meaning of $n$. From (3.10) and (3.11) we know that $n$ denotes the point in which the solution is evaluated, but when $R_1(t)$ is defined it seems that $n$ is the degree of the interpolation. Perhaps the authors are trying to build polynomials of increasing degree? The answer is negative because the construction of the method is different. It is clear that the authors make a great confusion between the mesh point at which the solution is approximated and the degree of the interpolant polynomial; the formulation of the remainder $R_1(t)$ is therefore illogical and completely wrong.
	
	\item In the middle of page 6 the third derivative of $f$ is assumed bounded (although in the statement of the Theorem it is just $f$ to be assumed bounded and not its derivatives). First of all, this is a quite unrealistic assumption since, as any researcher in this field should know, the solution of a FDE (and consequently the rhs of the FDE) does in general not have such kind of smoothness (e.g. see \cite{Diethelm2007,Lubich1983,Stynes2016}). However, and this is more important, why the authors consider the third derivative? Maybe they assume now an interpolant polynomial of degree 2? Clearly this is unreasonable since in the construction of the method in Eq. (3.6) linear polynomials are used. The confusion is even greater when derivatives of $f$ with order $3$ and order $n+1$ are used in the same equation to find a bound for $\|R_n(t)\|_{\infty}$ !
	
	\item The first bound at page 7 is completely illogical and wrong: it seems that each of the three terms $t-t_i$ (but actually they should be just two!) are bounded by $h$; since $t$ varies in $[0,t_{n+1}]$ they should be instead bounded by something roughly proportional to $t_{n+1}$ or $(n+1)h$.
	
	\item The correction of the mistakes listed in the previous point would lead to an error which is not bounded by something proportional to $h^{3+\alpha}n^{\alpha}$, as stated in Theorem 3.1, but (even when the upper bound 2 in the product in the last equations at page 6 is corrected to 1 as it should be) by something proportional to $t_n^{2+\alpha}$: this means that the method diverges, i.e. the error increases with $t_n$ and the method is completely unreliable and useless.
	
	\item From the middle of page 7 it seems that the proof is repeated by using different arguments in order to obtain different results: I do not comment on this further proof but I just point out that is not possible to understand why this proof is proposed since no corresponding statement is given in the theorem (a further sign of a very poor way of writing papers).
	
	\item The error analysis neglects the effect of the accumulation of the approximation errors: indeed, when the method is used in practice, one cannot use the previous exact value $y(t_n)$ (as stated in the formulation of the method) but the previous approximation $y_n$. 
	
\end{enumerate}

It is not worthwhile to study in detail the error of the scheme proposed in \cite{AtanganaOwolabi2018} and provide a correct and detailed analysis; it is however clear from the above reasoning, and this is the most general and important result, that the error is proportional to $t_n^{2+\alpha}$ and, therefore, not only the approximated solution does not converge to the exact solution as $h \to 0$ (as it should be required by any numerical method) but the error increases in a fast way as the integration moves forward: actually, the method diverges from the solution.

Methods devised according to the approach proposed in \cite[Eq. (3.11)]{AtanganaOwolabi2018} are therefore completely useless and dangerous and it is important to warn other scientists and researchers to not use them for their simulations: only wrong results can be indeed expected.

It is also worthwhile to point out that the paper \cite{AtanganaOwolabi2018} is written in a such poor way which can cause further confusion to readers. For instance, it is necessary to observe that:
\begin{itemize}
	\item the Definition 2.1 of the Caputo fractional derivative is written in a wrong way; it seems that the authors make confusion between the fractional derivative of $y(t)$ and the Riemann-Liouville integral of $f(t,y(t))$;
	\item although the authors claim to consider a fractional derivative of order $\alpha > 0$, from the initial conditions it seems that they actually restrict their investigation to the case $0<\alpha<1$;
	\item in the formulation of the method given in Eq. (3.11) the numerical approximation is denoted by the same notation $y(t_{n+1})$ used for the exact solution (instead of using a different notation such as $y_{n+1}$), thus misleading the readers which could be induced to believe that the scheme evaluates the exact solution and not its approximation; 
	\item in the notation $R^{\alpha}_n(t)$ used for the error in Theorem 3.1 it is reported a dependence on $t$ which actually does not exist;
	\item as a consequence of the (wrong) proof of Theorem 3.1 proposed by the authors, in the expression of the error $R^{\alpha}_n(t)$ the term $n^2$ should be replaced by $n^\alpha$;
	\item in the first (unnumbered) equation at page 4 and in the subsequent equation (3.5) a minus sign should be used instead of a plus sign; 
	\item the correct evaluation of $A_{\alpha,2}$ is Eq. (3.10) is
	\[
		A_{\alpha,2} = \frac{f(t_n,y_n)}{h\Gamma(\alpha)} \left\{ \frac{ht_n^{\alpha}}{\alpha} - \frac{t_n^{\alpha+1}}{\alpha+1}\right\} + \frac{f(t_{n-1},y_{n-1})}{h\Gamma(\alpha)(\alpha+1)} t_{n}^{\alpha+1}
	\]
	\item as a consequence of the above two points the formal derivation of the method is not the one presented in Eq. (3.11) but the one presented in this paper; 
	\item the rhs of the FDE is approximated by exactly the same interpolant polynomial in the two terms of (3.5) but the authors use $L_1(t)$ and $L_2(t)$ to denote them in (3.12) and the same different notation is used for the remainder; in particular, the same remainder is denoted in different parts of the proof as $R_1(t)$, $R_2(t)$ and $R(t)$; 
\end{itemize}

\section{A simple numerical experiment}

To provide a numerical evidence about the wrong approximation provided by the method described in \cite{AtanganaOwolabi2018} we consider the simple linear test problem
\[
	\DC^{\alpha}_0 y(t) = \lambda y(t) , \quad y(0) = y_0, \quad 0 < \alpha < 1
\]
whose exact solution can be expressed in terms of the Mittag-Leffler function according to 
\[
	y(t) = E_{\alpha}(t^{\alpha} \lambda) y_0
	, \quad E_{\alpha}(z) = \sum_{k=0}^{\infty} \frac{z^k}{\Gamma(\alpha k + 1)} .
\]

\begin{figure}[ht]
\centering \includegraphics[width=0.6\textwidth]{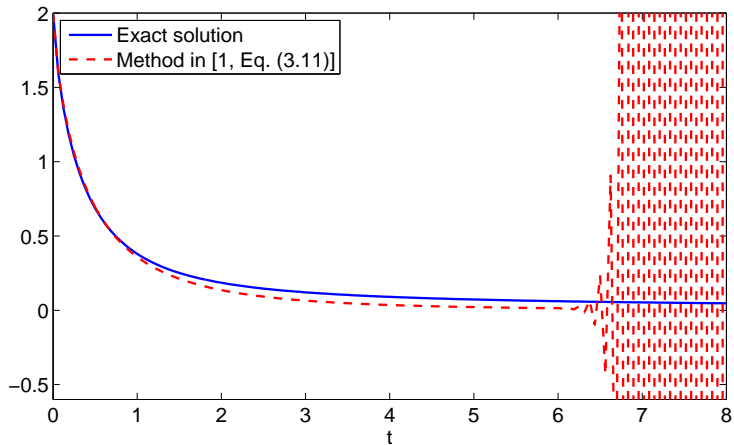}
\caption{Comparison of solutions for the linear test problem}\label{fig:Fig_AtanganaOwolabi_Comparison}
\end{figure}

In Figure \ref{fig:Fig_AtanganaOwolabi_Comparison} we compare, for $\alpha=0.8$, $\lambda = -2$ and $y_0=2$, the result obtained by the method proposed in \cite{AtanganaOwolabi2018} and the exact solution evaluated by means of the code devised in \cite{Garrappa2015_SIAM}. As we can clearly observe, the numerical solution diverges in a fast way from the exact solution, thus providing a meaningless approximation. 

\begin{figure}[ht]
\centering \includegraphics[width=0.6\textwidth]{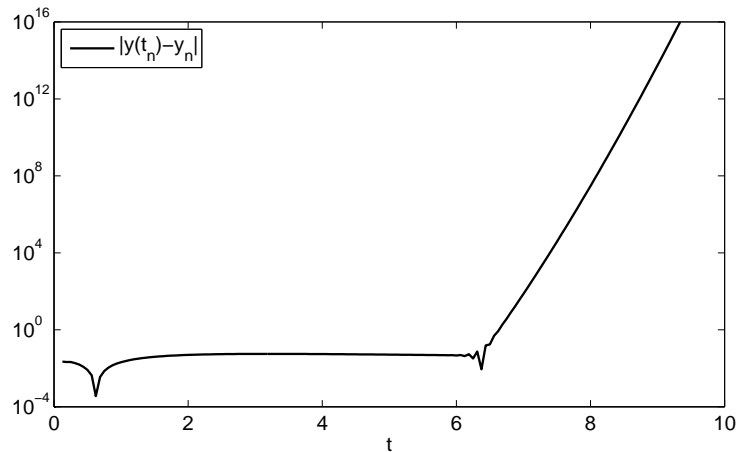}
\caption{Difference between exact and approximated solutions for the linear test problem}\label{fig:Fig_AtanganaOwolabi_Error}
\end{figure}

The plot of the error inf Figure \ref{fig:Fig_AtanganaOwolabi_Error} is even more clear in showing the divergence of the approximation from the exact solution (for these experiments we have used the formulation of the method presented in this paper since the original one published in \cite{AtanganaOwolabi2018} is affected by several misprints and errors and provides even wronger results).

These numerical findings have been obtained by selecting a step-size $h=2^{-4}$. Since the weak way by which the method has been derived, reducing the step-size does not help to obtain better results but just to increase the instability and divergence of the method.

\section{Further considerations}

The same approach used in \cite{AtanganaOwolabi2018} to derive the method for FDEs with the Caputo fractional derivative is also used, in the same paper, to derive methods for approximating the solution of equations with the so called Caputo-Fabrizio and ABC derivatives; in view of the marginal interest of these operators (e.g., see \cite{Giusti2018}) we do not discuss these methods but we just inform the readers that similar inconsistencies are encountered also in their derivation.


Unfortunately, it seems that the correct framework for the numerical treatment of fractional-order problems is not still widespread among non specialists; it is therefore possible that inappropriate, flawed or wrong methods, as the one proposed in \cite{AtanganaOwolabi2018}, are published in the literature and used for simulations. With this note we would recommend researchers to always check the reliability of the results presented in the literature, especially when they appear too much pretentious with respect to more established results, and encourage editors and reviewers to make a further effort in order to identify incorrect ideas and avoid their spread.

\end{document}